\documentclass[a4paper,12pt]{article}
\usepackage{amsmath,amsthm,amssymb,latexsym}
\title{{\bf On the semiadditivity of the capacities associated with signed vector valued Riesz kernels}}
\author{\Large{\Large Laura Prat}}
\date{}
\begin{document}
\maketitle
\newtheorem{teo}{Theorem}
\newtheorem{co}[teo]{Corollary}
\newtheorem{lemma}[teo]{Lemma}
\newtheorem{defi}[teo]{Definition}
\newtheorem{note}[teo]{Note}
\newtheorem{prop}[teo]{Proposition}
\newcommand{\Ha}{{\cal H}^{\alpha}}
\newcommand{\Hu}{{\cal H}^1}
\newcommand{\Rn}{{\mathbb R}^n}
\newcommand{\Rd}{{\mathbb R}^2}
\newcommand{\ep}{\varepsilon}
\newcommand{\N}{\mathbb{N}}
\newcommand{\Z}{\mathbb{Z}}
\newcommand{\C}{\mathbb{C}}
\newcommand{\cc}{{\cal C}}
\newcommand{\ga}{\gamma_\alpha}
\newcommand{\gai}{\gamma_{\alpha}^i}
\newcommand{\gam}{\gamma_{\alpha,+}}
\newcommand{\gaop}{\gamma_{\alpha,op}}
\newcommand{\al}{\alpha}
\begin{abstract}
 The aim of this paper is to show the semiadditivity of the
capacities  associated to the signed vector valued Riesz kernels of homogeneity
$-\al$ in $\Rn$, $0<\al<n$.
\end{abstract}

\section{Introduction}

%In this paper we show that, in $\Rn$, the capacity $\ga$ associated
%to the $\al-$dimensional Riesz kernel $x/|x|^{1+\al}$, $0<\al<n$, is
%semiadditive.
In this paper we study the capacity $\ga$ related to the signed
vector valued Riesz kernels $x/|x|^{1+\al}$  in $\Rn$, $0<\al<n$. If
$E\subset\Rn$ is a compact set and $0<\al<n$, one sets
\begin{equation}\label{ga}
\ga(E)=\sup|\langle T,1\rangle|,
\end{equation}
where the supremum is taken over all real distributions $T$ supported on
$E$ such that $x_i/|x|^{1+\al}*T$ is a function in $L^\infty(\Rn)$
and $\|x_i/|x|^{1+\al}*T\|_\infty\le 1$, for $1\leq i\le n$.

When $n=2$ and $\al=1$, by the celebrated result of X. Tolsa
\cite{semiad}, $\gamma_1$ is basically analytic capacity. Recall
that the analytic capacity of a compact set $E\subset\C$ is defined
as
\begin{equation}\label{acapacity}
\gamma(E)=\sup|\langle T, 1 \rangle|,
\end{equation}
the supremum taken over all complex distributions $T$ supported on
$E$ whose Cauchy potential $1/z*T$ is a bounded function and
$\|1/z*T\|_\infty\leq 1$.

The case $\al=n-1$, $n\ge 2$, is also particularly relevant, because
$\gamma_{n-1}$ coincides with the Lipschitz harmonic capacity, introduced in \cite{paramonov} to study problems of ${\cal
C}^1$-approximation by harmonic functions in $\Rn$ (see also
\cite{mp} and \cite{verderacm}).  Notice that the fact that,
in the plane, analytic capacity and $\gamma_1$ (Lipschitz harmonic
capacity) are comparable cannot be deduced just by an inspection of
(\ref{ga}) and (\ref{acapacity}). The reason is that the
distributions involved in the supremum in (\ref{acapacity}) are
complex.\newline

In \cite{laura1} one discovered the fact that if $0<\al<1$, then a
compact set of finite $\al$-dimensional Hausdorff measure has zero
$\ga$- capacity. This is in strong contrast with the situation for
integer $\al$, in which $\al$-dimensional smooth hypersurfaces have
positive $\ga$ capacity (see \cite{mp}). The case of non-integer $\al > 1$ is not
well understood, although it was shown in \cite{laura1} that for
Ahlfors-David regular sets the above mentioned result (for $0<\al<1$)
still holds in this case (see also \cite{laura3}).\newline

In \cite{mpv}, the surprising equivalence between $\ga$, $0 <\al < 1$, and one
of the well-known Riesz capacities of non-linear potential theory
(see \cite[ Chapter 1, p. 38]{adamshedberg}) was established. It was
shown that for some positive constant $C$,

\begin{equation}\label{crelle}
C^{-1}C_{\frac 2 3(n-\al),\frac 3 2}(E)\le \ga(E)\le CC_{\frac 2
3(n-\al),\frac 3 2}(E).
\end{equation}

Recall that the Riesz capacity $C_{s,p}$ of a compact set
$E\subset\Rn$, $1 < p < \infty$, $0 < sp\leq n$, is defined by

$$C_{s,p}(E)=\inf\{\|\varphi\|_p^p: \varphi*\frac 1{|x|^{n-s}}\geq 1\,\,\mbox{ on
}\,E\}$$ where the infimum is taken over all compactly supported
infinitely di�erentiable functions on $\Rn$. The capacity $C_{s,p}$
plays a central role in understanding the nature of Sobolev spaces
(see \cite{adamshedberg}).

In \cite{env} it has been shown that the first inequality in (\ref{crelle}) holds for all indices $0<\al<n$. The opposite inequality is false when $\al\in\Z$, for example if one takes $E$ contained in a $\al$-plane with positive $\al$-dimensional Hausdorff measure, then $\ga(E)>0$ while $C_{\frac 2 3(n-\al),\frac 3 2}(E)=0$. It is an open problem to prove (or
disprove) the second inequality in (\ref{crelle}) for non-integer $1<\al<n$.\newline

Since $C_{s,p}$ is a subadditive set function, as a direct consequence of (\ref{crelle}) one gets that $\ga$ is semiadditive, that is, given compact sets $E_1$
and $E_2$,

$$\ga(E_1\cup E_2)\leq C(\ga(E_1)+\ga(E_2)),$$
for some constant $C$ depending only on $n$ and $\al$. In fact $\ga$
is countably semiadditive. In this paper we will show that the semiadditivity of $\ga$ holds
for all indices $0<\al<n$ (see corollary $2$ below). \newline

If we restrict the supremum in (\ref{ga}) to distributions $T$ given
by positive Radon measures supported on $E$, we obtain the
capacities $\gam$. Clearly, we have
$$\gam(E)\leq\ga(E).$$

The arguments of this paper will prove that, in fact, these two
quantities are comparable, namely

\begin{teo}\label{alfa}
There exists an absolute constant $C>0$ such that for any compact
set $E\subset\Rn$ and any $0<\alpha<n$,
\begin{equation}\label{mainineq}
\ga(E)\leq C\gam(E).
\end{equation}
\end{teo}

%\begin{teo}\label{alfai}
%There exists an absolute constant $C>0$ such that for any compact
%set $E\subset\Rn$, any $0<\alpha\leq n-1$ and any component $1\leq
%i\leq n$, $$\gamma_\alpha^i(E)\leq C\gamma_{\alpha,+}^i(E).$$
%\end{teo}

This was first shown for $\al=1$ and $n=2$ by X. Tolsa
\cite{semiad}, and it was extended to the case $\al=n-1$ by Volberg
\cite{volberg}. For values of $\al\in(0,1)$ the result appears in
\cite{mpv}. In \cite{mt}, Theorem \ref{alfa} was proven for a
certain class of Cantor sets in $\Rn$, (see also \cite{tolsacantor}
where it is proven for a wider class of Cantor sets). Recently,
\cite{eidvol} have proven that this comparability result also holds on some
examples of random Cantor sets.\newline

As a corollary from Theorem \ref{alfa}, one deduces that $\ga$ is countably semiadditive for $0<\alpha<n$.\newline

\begin{co} \label{semiadit}
Let $E\subset\Rn$ be a compact set. Let $E_i$, $i\ge 1$, be Borel
sets such that $E = \bigcup_{i=1}^{\infty} Ei$. Then,
$$\ga(\bigcup_{i=1}^\infty E_i)\leq C\sum_{i=1}^\infty \ga(E_i)$$
where C is an absolute constant.
\end{co}

%To state our results in detail we need to introduce some notation.
%For $0<\al<n$, the $\al$-dimensional Riesz kernel is defined by
%$$K^\al(x) = \frac{x}{|x|^{\al+1}} ,\,\, x\in Rn,\,\, x\neq 0.$$

%\noindent Notice that this is a vectorial kernel. The
%$\al$-dimensional Riesz transform (or $\al-$Riesz transform) of a
%real Radon measure $\nu$ with compact support is

%$$R^\al\nu(x) = \int K^\al(y- x)d\nu(y),\, x\notin\mbox{supp}(\nu).$$

%\noindent Although the preceding integral converges a.e. with
%respect to Lebesgue measure, the convergence may fail for
%$x\in\,\mbox{supp}(\nu)$. This is the reason why one considers the
%truncated $\al$-Riesz transform of $\nu$ at level $\ep>0$, which is
%defined as

%$$R^\al_\ep\nu(x) = \int_{ |y-x|>\ep }K^\al(y-x) d\nu(y),\, x\in\Rn.$$

%\noindent These definitions also make sense if one considers
%distributions instead of measures. Given a compactly supported
%distribution $T$, set $$R^\al(T) = K^\al*T.$$

%\noindent (in the principal value sense for $\al = n$), and
%analogously $$R^\al_\ep(T) = K^\al_\ep*T,$$ where $K ^\al_\ep(x) =
%\chi_{|x|>\ep} x/|x|^{\al+1}$.

%Given a positive Radon measure $\mu$ with compact support and a
%function $f\in L^1(\mu)$, we consider the operators $R^\al_\mu(f):=
%R^\al(f d\mu)$ and $R^\al_{\mu,\ep}(f) := R^\al_\ep(f d\mu)$. We say
%that $R^\al_\mu$ is bounded on $L^2(\mu)$ when

%$$\|R^\al_\mu\|_{L^2(\mu)}=\sup_{\ep>0}\|R^\al_{\mu,\ep}\|_{L^2(\mu)}<\infty,$$
%or, in other words, when the truncated $\al$-Riesz transforms are
%uniformly bounded on $L^2(\mu)$.

The paper is organized as follows. In Section $2$ we prove Corollary
\ref{semiadit}. In Section $3$ we deal with one of the main
ingredients for the proof of Theorem \ref{alfa}, a localization
$L^\infty$-estimate for the scalar kernels $x_i/|x|^{1+\al}$,
$0<\al<n$, $1\leq i\leq n$. In Section $4$ we prove that the
capacities $\ga$ satisfy a exterior regularity property that will be
used for the proof of Theorem \ref{alfa}. Finally, in the last
section, we present a sketch of the proof of Theorem \ref{alfa}. It
becomes clear that the proof depends on the following three facts: a
localization $L^\infty$ estimate for the $\al-$Riesz kernels, the
exterior regularity property of $\ga$, $0<\al<n$, and Volberg's
extension \cite{volberg} of Tolsa's proof of the semiadditivity of
analytic capacity \cite{semiad}.\newline

Our notation and terminology are standard. For example, ${\cal
C}_0^\infty(E)$ denotes the set of all infinitely differentiable
functions with compact support contained in the set $E$. Cubes will
always be supposed to have sides parallel to the coordinate axis,
$l(Q)$ is the side length of the cube $Q$ and $|Q|=l(Q)^n$ its
volume.

Throughout all the paper, the letters $c,\,C$ will stand for
absolute constants depending only on $n$ and $\al$ that may
change at different occurrences.

\section{Proof of Corollary \ref{semiadit}.}

In this section we will deduce the semiadditivity of the $\ga$
capacity, $0<\alpha<n$, from Theorem \ref{alfa}. For this, we need
to introduce the $\al-$Riesz transform with respect to an underlying
positive Radon measure $\mu$ satisfying the $\al-$growth condition
\begin{equation}\label{mesgrowth}
\mu(B(x,r))\leq Cr^\al,\,\,x\in\Rn,\,\,r\ge 0.
\end{equation}

\noindent Given $\ep>0$ we define the truncated $\al-$Riesz transform at level $\ep$ as
$$R_\ep(f\mu)(x)=\int_{|y-x|>\ep}\frac{x-y}{|x-y|^{1+\al}}f(y)d\mu(y),\,\,x\in\Rn,$$
for $f\in L^2(\mu)$. The growth condition on $\mu$ insures that each $R_\ep$ is a bounded operator on $L^2(\mu)$ with operator norm $\|R_{\ep}\|_{L^2(\mu)}$ possibly depending on $\ep$. We say that the $\al-$Riesz transform is bounded on $L^2(\mu)$ when

$$\|R\|_{L^ 2(\mu)}=\sup_{\ep}\|R_\ep\|_{L^2(\mu)}<\infty,$$
or, in other words, when the truncated $\al-$Riesz transforms are uniformly bounded on $L^2(\mu)$. Call $L_\al(E)$ the set of positive Radon measures supported on $E$ which satisfy (\ref{mesgrowth}) with $C=1$. One defines $\gaop(E)$ by,

\begin{equation}\label{comppositiu}
\gaop(E)=\sup\{\mu(E):\,\mu\in L_\al(E)\,\mbox{ and }\,\|R\|_{L^2(\mu)}\leq 1\}.
\end{equation}

As it is well known,  the capacities $\gam(E)$ and $\gaop(E)$ are
comparable, that is, for some positive constant $C$ one has

\begin{equation}\label{gaopgam}
C^{-1}\gaop(E)\leq\gam(E)\leq C\gaop(E),
\end{equation}
for each compact set $E\subset\Rn$ (see Lemma 3 in \cite{laura3}).

Hence, once Theorem \ref{alfa} is available, namely the fact that $\ga(E)$ is comparable to $\gam(E)$, the semiadditivity of $\ga$ holds because $\gaop$ is
obviously semiadditive.

\section{Localization of $\al$-Riesz potentials}

\subsection{A growth condition}

 Let $T$ be a compactly supported distribution in $\Rn$ and
$0<\alpha<n$. Write $\al=[\al]+\{\al\}$, with $[\al]\in\Z$ and $0\le\{\al\}<1$. We say that the distribution $T$ has growth $\alpha$
provided
\begin{equation}\label{growthG}
G_\al(T) = \sup_{\varphi_Q} \frac{|\langle
T,\varphi_Q\rangle|}{l(Q)^{\al}} < \infty \,,
\end{equation}
where the supremum is taken over all $\varphi_Q \in
\cc^\infty_0(Q)$ satisfying the normalization inequalities

\begin{equation}\label{normal}
 \|\nabla^{n-\al}\varphi_Q\|\leq l(Q)^{\al},
\end{equation}
where $\|\nabla^{n-\al}\varphi_Q\|$ is defined as follows:

\begin{enumerate}

\item For $\al=[\al]\in\Z$, condition (\ref{normal}) means that
\begin{equation}\label{normalizationinteger}
\|\nabla^{n-\al}\varphi_Q\|:=\sup_{|s|=n-\al}\|\partial^s\varphi_Q\|_{L^1(Q)}
\leq C l(Q)^{\al}\,.
\end{equation}

\item for $\{\al\}>0$, condition (\ref{normal}) means that
\begin{equation}\label{normalization}
\|\nabla^{n-\al}\varphi_Q\|:=\sup_{|s|=n-[\al]}\|\partial^s\varphi_Q*\frac{1}{|x|^{n-\{\al\}}}\|_{L^1}
\leq C l(Q)^{\al}\,.
\end{equation}

\end{enumerate}

Here we are adopting the standard notation related to multi-indexes, that is, $s=(s_1,s_2,\cdots,s_n)$, where each coordinate $s_j$ is a non-negative integer and $|s|=s_1+\cdots+s_n$.

For a compact set $E$ in $\Rn$ we define $g_{\alpha}(E)$ as
the set of all distributions $T$ supported on $E$ having growth $\alpha$ with constant
$G_\alpha(T)$ at most $1$\,.\newline %For each coordinate $i$ set

%\begin{equation*}
%\gamma_{\alpha}^i(E)=\sup\{|<T,1>|\},
%\end{equation*}
%where the supremum is taken over those distributions $T\in
%g_\alpha(E)$, such that the $i$-th component of the $\alpha-$Riesz
%potential $f_i=x_i/|x|^{1+\alpha}*T$ is a bounded function and
%$\|f_i\|_\infty\leq 1$.

We start by showing that the usual $\al$-growth condition for a
positive Radon measure (see (\ref{mesgrowth})) is equivalent to the
notion of growth $\al$ for distributions, as defined in
(\ref{growthG}).

Given a positive Radon measure $\mu$ set

$$L_\al(\mu)=\sup_Q\frac{\mu(Q)}{l(Q)^\al}$$
where the supremum is taken over all cubes $Q$ with sides parallel
to the coordinate axis.

If $\varphi\in\cc^\infty_0(Q)$, then by an inequality of Maz'ya
\cite[p. 15 and p.134 ]{mazya}.

$$|\langle\mu,\varphi\rangle|=|\int\varphi d\mu|\leq\int|\varphi|d\mu\leq CL_\al(\mu)\|\nabla^{n-\alpha}\varphi\|.$$
Thus, $G_\al(\mu)\leq CL_\al(\mu).$

For the reverse inequality, given a cube $Q$, let $\varphi_Q$ be a
function in $\cc^\infty_0(2Q)$ such that $1\leq\varphi_Q$ on Q and
$\|\partial^s\varphi_Q\|_\infty\leq C_sl(Q)^{-|s|}$, $|s|\geq 0$.
Then  (\ref{normalizationinteger}) clearly holds when $\al\in\Z$. For $\{\al\}>0$ and $|s|=n-[\al]$, write $|r|=n-[\al]-1$. 
Bringing one derivative from $\partial^s\varphi_Q$ to the kernel to  get integrability in $(4Q)^c$ and using Fubini we obtain

\begin{equation*}
\begin{split}
\|\partial^s\varphi_Q*\frac 1{|x|^{n-\{\al\}}}\|_{L^1}&= \int_{(4Q)^c}\left(\partial^s\varphi_Q*\frac 1{|x|^{n-\{\al\}}}\right)(y)dy +
\int_{4Q}\left(\partial^s\varphi_Q*\frac 1{|x|^{n-\{\al\}}}\right)(y)dy\\\\& 
\leq \frac{C}{l(Q)^{|r|}}\int_{(4Q)^c}\int_Q\frac{dydx} {|y-x|^{n+1-\{\al\}}}+\frac{C}{l(Q)^{|s|}}\int_{4Q}\int_Q\frac{dydx}{|y-x|^{n-\{\al\}}}\\\\&
\leq Cl(Q)^{\alpha}.
\end{split}
\end{equation*}
Thus, $$\|\nabla^{n-\al}\varphi_Q\|:=\sup_{n-[\al]}\|\partial^s\varphi_Q*\frac 1{|x|^{n-\{\al\}}}\|_{L^1}\leq Cl(Q)^\al.$$

Therefore,

$$\mu(Q)\leq\int\varphi_Q d\mu\le |\langle\mu,\varphi_Q\rangle|\leq C\,G_\al(\mu)\,l(Q)^{\al}.$$

Next lemma shows that all distributions admissible in the definition of $\ga(E)$ have growth $\al$.

\begin{lemma}\label{growthlemma}
Let $T$ be a distribution supported on the compact set
$E\subset\Rn$. Let $0<\alpha<n$ and suppose that $T$ has bounded
$\alpha$-Riesz potential $x/|x|^{1+\al}*T$. Then $T\in g_\al(E)$.
\end{lemma}

{\em Proof.} Our proof uses a reproduction formula for test functions involving the kernel $x_i/|x|^{1+\al},$ which was first introduced in 
\cite[Lemma 3.1]{laura1}. There are many variants of this formula depending, for instance, on whether de dimension $n$ and the integer part of $\al$ are 
even or odd. We will consider in full detail only the case of odd dimension of the form $n=2k+1$. We will also assume that $\al$ is non-integer and that 
its integer part is even, of the form $[\al]=2d$. The argument for the other cases follows the same line of reasoning but using the different variants of the 
corresponding reproduction formula. In our present case, the reproduction formula we use reads as follows,

\begin{equation}\label{betterrepro}
\varphi(x)=c\sum_{i=1}^n\Delta^{k-d}\partial^i\varphi*\frac{x_i}{|x|^{1+\al}}*\frac 1{|x|^{n-\{\al\}}}.
\end{equation}

Let $\varphi_Q$ be a ${\cal C}_0^{\infty}(Q)$ function satisfying
the normalization inequalities (\ref{normal}). Then, by
(\ref{betterrepro}), the boundedness of the potential
$x_i/|x|^{1+\al}*T$, $1\le i\le n$,  and Fubini

\begin{eqnarray*}
|\langle T,\varphi_Q\rangle|&\leq&
\,\sum_{i=1}^n|\langle \frac{x_i}{|x|^{1+\al}}*T,\Delta^{k-d}\partial_i\varphi_Q*\frac 1{|x|^{n-\{\al\}}}\rangle|\\\\
&\leq& C\sum_{i=1}^n\int|\Delta^{k-d}\partial_i\varphi_Q*\frac{1}{|x|^{n-\{\al\}}}(y)|dy\\\\&\leq&C\sum_{i=1}^n\int_{2Q}\int_Q\frac{|\Delta^{k-d}\partial_i\varphi_Q(z)|}{|z-y|^{n-\{\al\}}}dzdy +C\int_{(2Q)^c}\int_Q\frac{|\Delta^{k-d}\varphi_Q(z)|}{|z-y|^{n+1-\{\al\}}}dzdy\\\\&\leq& Cl(Q)^\alpha.
\end{eqnarray*}

The other cases (namely for odd $n$, even $[\al]$...)  are proven in
the same way by using analogous formulae (see \cite[Lemma
3.1]{laura1}).\qed

\subsection{Localization of Riesz potentials.}

When analyzing the argument for the proof of the semiadditivity of analytic capacity (see
Theorem 1.1 in \cite{semiad}) one realizes that one of the technical
tools used is the fact that the Cauchy kernel~$1/z$ localizes in the
uniform norm. By this we mean that if $T$ is a compactly supported
distribution such that $1/z*T$ is a bounded measurable function,
then $1/z*(\varphi \, T)$  is also bounded measurable for each
compactly supported ${\mathcal C}^1$ function~$\varphi$. This is an old
result, which is simple to prove because $1/z$ is related to the
differential operator~$\overline{\partial}$ (see~\cite[Chapter V]{garnett}). The same localization result can be proved easily in any
dimension for the kernel~$x/|x|^n$, which is, modulo a
multiplicative constant, the gradient of the fundamental solution of
the Laplacian. Again the proof is reasonably straightforward because
the kernel is related to a differential operator (see~\cite{paramonov} and~\cite{verderacm}).

In \cite[Lemma 3.1]{mpv} we  were concerned with the localization of
the vector valued $\alpha$-Riesz kernel $x/|x|^{1+\alpha}$,
$0<\alpha<n$. For general values of $\alpha$  there is no
differential operator in the background and consequently the
corresponding localization result becomes far from obvious. We state now the general localization Lemma proved in \cite{mpv}.\newline

In what follows, given a cube $Q$, $\varphi_Q$ will denote and infinitely differentiable function supported on $Q$ and such that $\|\partial^s\varphi_Q\|_\infty\le l(Q)^{-|s|}$, $0\le |s|\le n-[\al]$.\newline

\begin{lemma}\label{localizationmpv}
Let $T$ be a compactly supported distribution in $\Rn$ and let $0<\al<n$. Suppose that $x_i/|x|^{1+\alpha}*T$ is a bounded measurable function for $1\le i\le n$. Then there exists some constant $C=C(n,\al)>0$ such that
$$\sup_{1\le i\le n}\|\frac{x_i}{|x|^{1+\alpha}}*\varphi_Q T\|_\infty\leq C\sup_{1\le i\le n}\|\frac{x_i}{|x|^{1+\alpha}}* T\|_\infty.$$
\end{lemma}

Although Lemma \ref{localizationmpv} is enough for our purposes,
that is to prove Theorem \ref{alfa}, in this paper we will give a
proof of a stronger localization result, with a shorter and less
technically involved proof. The main difference between the
localization lemma in \cite{mpv} and the one we prove here is that
we localize one component of the vector potential
$\frac{x}{|x|^{1+\alpha}}*T$, only assuming $L^\infty$ estimates on
the potential of the same component, instead of assuming
$\|\frac{x}{|x|^{1+\alpha}}*T\|_\infty\leq 1$ for the whole vector.
Our new localization lemma reads as follows,

\begin{lemma}\label{localization1}
 Let $T$ be a compactly supported distribution in $\Rn$ with $\alpha$-growth, $0<\alpha<n$,
 such that $(x_i / |x|^{1+\alpha}) *T$ is in $L^\infty(\Rn)$ for some $i$\,, $1\leq i \leq n$.
 Then $(x_i / |x|^{1+\alpha}) * \varphi_Q T$ is in $L^\infty(\Rn)$ and
$$\|\frac{x_i}{|x|^{1+\alpha}}*\varphi_Q T\|_\infty\leq C\,(\| \frac{x_i}{|x|^2}*T\|_\infty+G_{\alpha}(T))\,,$$
for some positive constant $C=C(n)$ depending only on $n\,.$
\end{lemma}

For $\al=1$ the proof of the above lemma can be found in
\cite{mpv2}. We remark here that  when one deals with indexes $\alpha\in\Z$,  the proof of Lemma \ref{localization1} 
is less technically involved, since the derivatives $\partial^s\varphi_Q$, $|s|=n-\alpha$, are ordinary 
derivatives and therefore supported on the cube $Q$   (compare (\ref{normalization}) with 
(\ref{normalizationinteger})). \newline

For the proof of Lemma \ref{localization1} we need the following result (see Lemma $7$ in \cite{mpv2} for the case $\al=1$), that will 
be proved after the proof of Lemma \ref{localization1}.

\begin{lemma}\label{prelocalization}
 Let $T$ be a compactly supported distribution in $\Rn$ with $\alpha-$growth, $0<\alpha<n$.
 Then, for each coordinate $i$, the distribution $(x_i / |x|^{1+\alpha}) * \varphi_Q T$ is a locally integrable function
 and there exists a point $x_0 \in \frac{1}{4}Q$
  such that $$\left|\left(\frac{x_i}{|x|^{1+\alpha}}*\varphi_QT\right)(x_0)\right|\leq C \, G_{\alpha}(T)\,,$$
where  $C=C(n)$ is a positive constant depending only on $n\,.$
\end{lemma}

\vspace{.5cm}
\noindent {\em Proof of Lemma \ref{localization1}.} Without loss of
generality take $i=1$. %Since $k^1* \varphi_Q T$ is a harmonic
%function off the closure of $Q$, by the maximum principle we only
%need to estimate $|(k^1* \varphi_Q T) (x)|$ for $x\in \frac 3 2 Q$.
We distinguish two cases:
\begin{enumerate}
\item $x\in (\frac 3 2 Q)^c$. Then $ |(k^1*\varphi_QT)(x)|=|<T,\varphi_Q(y)k^1(x-y)>|.$
Notice that, for an appropiate dimensional constant $C$,  the
function $$\psi_Q(y)=Cl(Q)^\alpha\varphi_Q(y)k^1(x-y),$$ satisfies
the normalization inequalities in the definition of $G_\alpha(T)$,
namely
\begin{equation}\label{normaliz}
\|\nabla^{n-\al}\psi_Q\|\leq l(Q)^{\al}\,.
\end{equation}
Therefore,

$$|(k^1*\varphi_QT)(x)|=Cl(Q)^{-\alpha}|<T,\psi_Q>|\leq C.$$

To see (\ref{normaliz}), observe that if $[\al]$ denotes the integer
part of $\al$ and we write $\al=[\al]+\{\al\}$, then by Leibniz
formula,
\begin{equation}\label{psiq}
\|\partial^s\psi_Q\|_{L^\infty(Q)}\leq
Cl(Q)^\al\sum_{|r|=0}^{|s|}l(Q)^{-|r|}l(Q)^{-|s|-\al+|r|}\leq
Cl(Q)^{-|s|},
\end{equation}
for any multiindex $s=(s_1,\cdots,s_n)$ with $|s|\ge 0$.

If $\{\al\}=0$, (\ref{psiq}) immediately implies that condition
(\ref{normaliz}) holds. When
$\{\al\}>0$, let $s=(s_1,\cdots,s_n)$ be any multiindex with $|s|=n-[\al]$ and write

\begin{equation*}
\begin{split}
\int|(\partial^s\psi_Q*\frac 1{|z|^{n-\{\al\}}})(y)|dy&=\int_{2Q}|(\partial^s\psi_Q*\frac{1}{|z|^{n-\{\al\}}})(y)|dy\\\\&+
\int_{(2Q)^c}|(\partial^s\psi_Q*\frac{1}{|z|^{n-\{\al\}}})(y)|dy=A+B.
\end{split}
\end{equation*}
By (\ref{psiq}), we have
$$A\le C\|\partial^s\psi_Q\|_{\infty}\int_{2Q}\int_Q\frac{dzdy}{|z-y|^{n-\{\al\}}}\le Cl(Q)^\al.$$
And by bringing one derivative from $\partial^s\psi_Q$ to the kernel $|z|^{-n+\{\al\}}$ and using (\ref{psiq}) again, we get
\begin{equation}
 \begin{split}
  B&\le C\int_{(2Q)^c}|(\partial^t\psi_Q*\frac{1}{|z|^{n+1-\{\al\}}})(y)|dy\\\\&\le C\|\partial^t\psi_Q\|_{\infty}\int_{(2Q)^c}\int_Q\frac{dzdy}{|z-y|^{n+1-\{\al\}}}
\le Cl(Q)^\al,
 \end{split}
\end{equation}
for some multiindex $t$ with $|t|=n-[\al]-1$.

\item $x\in \frac 3 2 Q$. Since
$k^1*T$ and $\varphi_Q$ are bounded functions, we can write
\begin{eqnarray*}
 |(k^1*\varphi_QT)(x)|\leq|(k^1*\varphi_QT)(x)-\varphi_Q(x)(k^1*T)(x)|+\|\varphi_Q\|_\infty\|k^1*T\|_\infty.
\end{eqnarray*}

Let $\psi_Q\in{\cal C}_0^{\infty}(\Rn)$ be such that $\psi_Q\equiv
1$ in $2Q$, $\psi_Q\equiv 0$ in $(4Q)^c$ and
$\|\partial^s\psi_Q\|_\infty\leq C_s\,l(Q)^{-|s|}\,,$ for each
multi-index $s$ . Then one is tempted to write
\begin{eqnarray*}
 |(k^1*\varphi_QT)(x)-\varphi_Q(x)(k^1*T)(x)|&\leq&|<T,\psi_Q(y)(\varphi_Q(y)-\varphi_Q(x))k^1(x-y)>|\\\\&+&
\|\varphi_Q \|_{\infty}\,|<T,(1-\psi_Q(y))k^1(x-y)>|\,.
\end{eqnarray*}
The problem is that the first term in the right hand side above does
not make any sense because $T$ is acting on a function of $y$ which
is not necessarily differentiable at the point $x\,.$  To overcome
this difficulty one needs to use a standard regularization process.
Take $\chi \in \cc^\infty(B(0,1))$ such that $\int \chi(x)\,dx = 1$
and set $\chi_\ep(x)= \ep^{-n}\,\chi(x/\ep)\,.$ The plan is to
estimate, uniformly on $x$ and $\epsilon\,,$
\begin{equation}\label{reg}
|(\chi_\ep*k^1*\varphi_QT)(x)-\varphi_Q(x)(\chi_\ep*k^1*T)(x)|\,.
\end{equation}
 Clearly \eqref{reg} tends, as $\ep$ tends to zero,
 to
$$
|(k^1*\varphi_QT)(x)-\varphi_Q(x)(k^1*T)(x)|\,,
$$
for almost all $x \in \Rn$\,, which allows the transfer of uniform
estimates. We now have
\begin{eqnarray*}\label{dif}
|(\chi_\ep*k^1*\varphi_QT)(x)&-&\varphi_Q(x)(\chi_\ep*k^1*T)(x)|\\\\&\le&
|<T,\psi_Q(y)(\varphi_Q(y)-\varphi_Q(x))
(\chi_\ep*k^1)(x-y)>|\\\\&+&
\|\varphi_Q\|_{\infty}|<T,(1-\psi_Q(y))(\chi_\ep*k^1)(x-y)>|\\\\&=&A_1+A_2.
\end{eqnarray*}
\noindent where the last identity is the definition of $A_1$ and
$A_2$. \noindent To deal with term $A_1$ set
$$k^{1,x}_\ep(y)=(\chi_\ep*k^1)(x-y).$$ We claim that, for an
appropriate dimensional constant $C$, the test function
$$f(y)=Cl(Q)^\al\psi_Q(y)(\varphi_Q(y)-\varphi_Q(x))k_\ep^{1,x}(y),$$
satisfies the normalization inequalities \eqref{normal} in
the definition of $G_\al(T)\,,$ with $\varphi_Q$ replaced by $f$. If
this is the case, then
$$A_1\leq C l(Q)^{-\al}|<T,f>|\leq C\,G_\al(T).$$

To prove the claim we have to show that
\begin{equation}\label{regdins}
\|\nabla^{n-\al}f\|\leq
l(Q)^{\al}.
\end{equation}

We first notice that the regularized kernel $\chi_\ep*k^1$ satisfies
the inequalities
\begin{equation}\label{regkernel}
|(\chi_\ep* \partial^s \,k^1)(x)| \le \frac{C}{|x|^{\al+|s|}}\,,
\quad x \in \Rn\setminus \{0\}\quad \text{and}\quad 0\le |s| \leq
n-[\alpha]-1\,,
\end{equation}
where $C$ is a dimensional constant, which, in particular, is
independent of $\epsilon$.

The estimate of the $L^1-$norm in (\ref{regdins}) requires the use
of Leibniz formula
\begin{equation}\label{Leibniz2}
\partial^s \left(\psi_Q(\varphi_Q -\varphi_Q(x))k_{\ep}^{1,x}\right) =
\sum_{|r|=0}^{|s|}c_{r,s}\,
\partial^r(\psi_Q(\varphi_Q
-\varphi_Q(x)))\;\partial^{s-r}\, k_\ep^{1,x}
\end{equation}
and of \eqref{regkernel}\,. Notice that for any multiindex $t=(t_1,\cdots,t_n)$ with $|t|=n-[\al]-1$,
\begin{equation}\label{partials}
\begin{split}
\|\partial^t f \|_{L^1(4Q)}&\le C l(Q)^\al\sum_{|r|=0}^{|t|}\frac
1{l(Q)^{|r|}}\int_{4Q}|\partial^{t-r}(k_\ep^{1,x})(y)|\,dy \\\\
&\leq C l(Q)^{n-|t|}=Cl(Q)^{[\al]+1}.
\end{split}
\end{equation}

And if $s=(s_1,\cdots, s_n)$ is such that $|s|=n-[\al]$, using the mean value theorem to gain integrability when $|r|=0$ and (\ref{regkernel}),
\begin{equation}\label{partialt}
\begin{split}
 \|\partial^s f \|_{L^1(4Q)}&\le Cl(Q)^{\al}\|\nabla\varphi_Q\|_\infty\int_{4Q}\frac{dy}{|y-x|^{\al+|s|-1}} \\\\&+C l(Q)^\al\sum_{|r|=1}^{|s|}\frac
1{l(Q)^{|r|}}\int_{4Q}|\partial^{s-r}(k_\ep^{1,x})(y)|\,dy\\\\&\leq
Cl(Q)^{[\al]}.
\end{split}
\end{equation}

Estimate (\ref{partialt}), immediately yields
(\ref{regdins}) for $\{\al\}=0$,
$$\|\nabla^{n-\al}f\|:=\sup_{|s|=n-\al}\int_{4Q}|\partial^sf(y)|dy\leq
Cl(Q)^{\al}.$$

If $\{\al\}>0$, then for any multiindex $s=(s_1,\cdots,s_n)$ with $|s|=n-[\al],$
\begin{eqnarray*}
\int|\,(\partial^s f*\frac{1}{|z|^{n-\{\al\}}})(y)|\,dy&\le&
\int_{5Q}|\,(\partial^s f*\frac{1}{|z|^{n-\{\al\}}})(y)|\,dy\\\\&+
&C\int_{(5Q)^c}|\,(\partial^t f*\frac{1}{|z|^{n+1-\{\al\}}})(y)|\,dy\\\\&=&B_1+B_2,
\end{eqnarray*}
where $t$ is some multiindex with $|t|=n-[\al]-1$.

To estimate $B_2$, we use Fubini and (\ref{partials}).  Then,
$$B_2\leq\int_{4Q}|\partial^tf(z)|\int_{(5Q)^c}\frac
{dydz}{|z-y|^{n+1-\{\al\}}}\leq Cl(Q)^{\al}.$$

We turn now to term $B_1$. By Fubini and (\ref{partialt}), we get
$$B_1\leq C\int_{5Q}|\partial^sf(z)|\int_{4Q}\frac{dydz}{|y-z|^{n-\{\al\}}}\leq
Cl(Q)^{\al}.$$

This finishes the proof of (\ref{regdins}). We now turn to $A_2$. By
Lemma \ref{prelocalization}, there exists a point $x_0\in Q$ such
that $|(k^1*\psi_QT)(x_0)|\leq C\, G_\al(T)$.
 Then
$$|(k^1*(1-\psi_Q)T)(x_0)|\leq C\,(\|k^1*T\|_\infty +G_\al(T)).$$
The analogous inequality holds as well for the regularized
potentials appearing in $A_2\,,$ uniformly in $\epsilon\,,$ and
therefore

$$A_2\leq C\,|<T,(1-\psi_Q)(k_\ep^{1,x}-k_\ep^{1,x_0})|+C\,(\|k^1*T\|_\infty +G_\al(T)).$$

To estimate $|<T,(1-\psi_Q)(k_\ep^{1,x}-k_\ep^{1,x_0})>|$, we
decompose $\Rn \setminus \{x\}$ into a union of rings $$R_j=\{z\in
\Rn:\;2^j\,l(Q)\leq|z-x|\leq 2^{j+1}\,l(Q)\},\;\; j\in\mathbb{Z},$$
and consider functions $\varphi_j$ in ${\cal C}^\infty_0(\Rn)$, with
support contained in $\frac 3 2 R_j$\,, such that
$\|\partial^s\varphi_j\|_\infty\leq C \,(2^j\,l(Q))^{-|s|}\,,$ $|s|
\geq 0$, and $\sum_j\varphi_j=1$ on $\Rn\setminus\{x\}$. Then, since
$x\in\frac 3 2 Q$ and $1-\psi_Q\equiv 0$ in $2Q$, the smallest ring
$R_j$ that may intersect $(2Q)^c$ is $R_{-2}$. Therefore  we have

\begin{eqnarray*}
 |<T,(1-\psi_Q)(k_\ep^{1,x}-k_\ep^{1,x_0})>|&=&<T,\sum_{j\geq -2}\varphi_j(1-\psi_Q)(k_\ep^{1,x}-k_\ep^{1,x_0})>|\\\\
 &\leq&|<T,\sum_{j\in I}\varphi_{j}(1-\psi_Q)(k_\ep^{1,x}-k_\ep^{1,x_0})>|\\\\&+&\sum_{j\in J}|<T,\varphi_{j}(k_\ep^{1,x}-k_\ep^{1,x_0})>|,
\end{eqnarray*}
where $I$ denotes the set of indices $j\geq -2$ such that the
support of $\varphi_j$ intersects $4Q$  and $J$ the remaining
indices, namely those $j \geq -2 $ such that the support of
$\varphi_j$ is contained in the complement of $4Q\,.$ Notice that
the cardinality of $I$ is bounded by a dimensional constant.
\newline

\noindent Set
$$g =C\,l(Q)^\al\,\sum_{j\in I}\varphi_j(1-\psi_Q)\,(k_\ep^{1,x}-k_\ep^{1,x_0})\,,$$
and for $j\in J$

$$g_j=C\,2^j\,(2^jl(Q))^\al\,\varphi_j\,(k_\ep^{1,x}-k_\ep^{1,x_0}).$$
We show now that the test functions $g$ and $g_j$, $j\in J\,,$
satisfy the normalization inequalities \eqref{normal} in the
definition of $G_\al(T)\,.$ Once this is available, using the $\alpha$-
growth condition of $T$ we obtain
\begin{eqnarray*}
 |<T,(1-\psi_Q)(k_\ep^{1,x}-k_\ep^{1,x_0})>|&\leq& C l(Q)^{-\al}|<T,g>|\\\\&+&
 C \sum_{j\in J} 2^{-j}(2^jl(Q))^{-\al}|<T,g_j>|\\\\
 &\leq& C\,G_\al(T) + C\sum_{j\geq 2}2^{-j}\,G_\al(T)\leq C\,G_\al(T)\,,
\end{eqnarray*}
which completes the proof of Lemma 7.

We check now the normalization inequalities for $g$ and $g_j$.
For $g$ one argues as in the proof of (\ref{regdins}), using that $\|\partial^s(1-\psi_Q)\|_\infty\leq
Cl(Q)^{-|s|},$ $\|\partial^s\varphi_j\|_\infty\leq C\,
l(Q)^{-|s|}$, $j\in I$,  \eqref{regkernel}, the
fact that  $x, x_0 \in \frac 3 2 Q$\,,\,$y\in (2Q)^c$\,, and a
gradient estimate.

For $g_j$ we use in addition Leibniz formula and
a gradient estimate to show that, for $j\in J$, and $n-[\al]-1\le|s|\le n-[\al]$,

\begin{equation}\label{gj}
\begin{split}
\|\partial^s g_j\|_{\infty}&\leq C\,
2^j(2^{j}\,l(Q))^\al\sum_{|r|=0}^{|s|}\frac
1{(2^jl(Q))^{|r|}}\frac{l(Q)}{(2^j\,l(Q))^{1+\al+|s|-|r|}}
\\&\leq C\,(2^jl(Q))^{-|s|},
\end{split}
\end{equation}

\noindent If $\{\al\}=0$,  we use (\ref{gj}) to obtain
$$\|\nabla^{n-\al}g_j\|:=\sup_{|s|=n-\al}\int|\partial^s g_j(y)|dy\le C(2^jl(Q))^{\al}.$$

\noindent If $\{\al\}>0$, then for any multiindex $s=(s_1,\cdots,s_n))$ with $|s|=n-[\al]-1$, 
using Fubini, (\ref{gj}) and arguing similar to the proof of (\ref{regdins}) we get,
for some multiindex $t=(t_1,\cdots, t_n)$ with $|t|=n-[\al]-1$,

\begin{equation*}
\begin{split}
\int|(\partial^sg_j*\frac 1{|z|^{n-\{\al\}}})(y)|dy&\le\int_{2R_j}|(\partial^sg_j*\frac 1{|z|^{n-\{\al\}}})(y)|dy\\\\&+
C\int_{(2R_j)^c}|(\partial^t g_j*\frac 1{|z|^{n+1-\{\al\}}})(y)|dy\\\\&\le C(2^jl(Q))^\al.
\end{split}
\end{equation*}

Therefore, we can conclude that $$\|\nabla^{n-\al}g_j\| \le C(2^jl(Q))^\al.\qed$$

\end{enumerate}

{\em Proof of Lemma \ref{prelocalization}.} Without loss of generality set $i=1$ and write
$k^1(x)=x_1/|x|^{1+\al}$. Since $k^1 * \varphi_Q T$ is infinitely
differentiable off the closure of ${Q}\,,$ we only need to show that
$k^1* \varphi_Q T$ is integrable on $2Q\,.$ We will actually prove a
stronger statement, namely, that $k^1* \varphi_Q T$ is in $L^p(2Q)$
for each $p$ in the interval $1\leq p< n\,.$ Indeed, fix any $q$
satisfying $n/(n-1) <q < \infty$ and call $p$ the dual exponent, so
that $1 < p< n\,.$ We need to estimate the action of $k^1 *
\varphi_Q T$ on functions $\psi \in \cc^\infty_0(2Q)$ in terms of
$\|\psi\|_q \,.$ We clearly have
$$
< k^1 * \varphi_Q T, \psi> = <T, \varphi_Q\,(k^1 * \psi)> \,.
$$
 We claim that, for an appropriate
dimensional constant $C \,,$ the test function
$$\frac{\varphi_Q\,(k^1 * \psi)}{C \,l(Q)^{\frac n p-\al} \,\|\psi\|_q}   $$
satisfies the normalization inequalities \eqref{normal} in
the definition of $G_\alpha(T)\,.$  Once this is proved, by the
definition of $G_\alpha(T)$ we get
\begin{equation*}\label{Lq}
|< k^1 * \varphi_Q T, \psi>| \le C\, l(Q)^\frac n p\,\|\psi\|_q
\,G_\alpha(T)\,,
\end{equation*}
and so
\begin{equation*}\label{Lp}
\|k^1 * \varphi_Q T \|_{L^p(2Q)} \le C\, l(Q)^\frac n
p\,G_\alpha(T)\,.
\end{equation*}
Hence
\begin{equation*}
\begin{split}
\frac{1}{|\frac{1}{4}Q|}\,\int_{\frac{1}{4} Q} |(k^1 * \varphi_Q
T)(x)|\,dx &\le 4^n\,\left(\frac{1}{|Q|}\,\int_Q |(k^1 * \varphi_Q
T)(x)|^p\,dx\right)^{\frac 1 p}\\& \le
C\,G_\alpha(T)\,,
\end{split}
\end{equation*}
which completes the proof of Lemma 8.

To prove the claim we have to show that
\begin{equation}\label{condition}
\| \nabla^{n-\al}(\varphi_Q\,(k^1 * \psi))\| \le C\,
l(Q)^{\frac n p}\, \|\psi\|_q.
\end{equation}

Write $\al=[\al]+\{\al\}$, with
$\{\al\}\in[0,1)$ and $[\al]\in\Z$. We distinguish now two cases, $\{\al\}$=0 and $\{\al\}>0$. 

\begin{enumerate}
\item Case $\{\al\}$=0, i.e. $\al=[\al]\in\Z$. This is the easiest case, because the derivatives appearing in \eqref{condition} are 
ordinary derivatives (see also lemma $7$ in \cite{mpv2}).

Let $s=(s_1,s_2,\cdots,s_n)$ be any multiindex with $|s|=n-\al$. Using Leibniz formula, 
\begin{equation}\label{Leibniz}
\partial^s \left(\varphi_Q \,(k^1 * \psi)\right) =  \sum_{|r|=0}^{|s|}
c_{s,r}\,\partial^r\varphi_Q \;\partial^{s-r}(k^1*\psi)\,,
\end{equation}
we obtain

\begin{equation*}
 \begin{split}
  \int_{Q}|\partial^s(\varphi_Q(k^1*\psi))(y)|dy&\le
C\int_{2Q}|(\partial^s k^1*\psi)(y)|dy\\\\&+C\sum_{|r|=1}^{|s|}\int_{Q}|\partial^r\varphi_Q(y)||\partial^{s-r}(k^1*\psi)(y)|dy=A+B.
 \end{split}
\end{equation*}

To estimate term $A$, we remark that for
$|s|=n-\alpha$,
$$\partial^s k^1*\psi=c\psi+T(\psi),$$ where $T$ is a smooth
homogeneous convolution Calder\'on-Zygmund operator and $c$ a
constant depending on $s$. This can be seen by computing the Fourier
transform of $\partial^s k^1$ and then using that each homogeneous
polynomial can be decomposed in terms of homogeneous harmonic
polynomials of lower degrees ( see \cite[p. 69]{St}). Since
Calder\'on-Zygmund operators preserve $L^q(\Rn)$, $1<q<\infty$, we
get, using H\"older
$$A\le Cl(Q)^{\frac n p}\|\psi\|_q.$$

To estimate term $B$, we use $|\partial^{s-r}k^1(x)| \le C\, |x|^{-(\alpha+|s|-|r|)}\,,$ Fubini, the fact that $\|\partial^r\varphi_Q\|_\infty\leq C
l(Q)^{-|r|}$ and H\"older to obtain $B\le Cl(Q)^{\frac n p}\|\psi\|_q$. Therefore we get,
$$\| \nabla^{n-\al}(\varphi_Q\,(k^1 * \psi))\|=\sup_{|s|=n-\al}\int_{Q}|\partial^s(\varphi_Q(k^1*\psi))(y)|dy\le C l(Q)^{\frac n p}\|\psi\|_q. $$

\item Case $\{\al\}>0$. Let $s=(s_1,s_2,\cdots,s_n)$ be any 
multiindex with $|s|=n-\al$ and write

%\begin{eqnarray*}
%\int|\nabla^{n-\alpha}(\varphi_Q(k^1*\psi))(x)|dx&=&\int_{(2Q)^c}|\partial^{s}(\varphi_Q(k^1*\psi))*\frac
%1{|x|^{n+1-\{\al\}}})(x)|dx\\\\&+&\int_{2Q}|\partial^{s-1}(\varphi_Q(k^1*\psi))*\frac
%1{|x|^{n-\{\al\}}})(x)|dx\\\\&=&A+B,
%\end{eqnarray*}
%the last identity being a definition for $A$ and $B$.

\begin{equation*}
\begin{split}
 \int|\partial^s(\varphi_Q(k^1*\psi)*\frac 1{|x|^{n-\{\al\}}})(x)|dx&=
 \int_{(2Q)^c}|(\partial^s \varphi_Q(k^1*\psi)*\frac 1{|x|^{n-\{\al\}}})(x)|dx\\\\&+
 \int_{2Q}|\partial^s(\varphi_Q(k^1*\psi)*\frac 1{|x|^{n-\{\al\}}})(x)|dx\\\\&=A+B.
\end{split}
\end{equation*}

\noindent We deal first with term $A$. Bringing one derivative from $\partial^s(\varphi_Q(k^1*\psi))$ to the kernel 
$|x|^{-n+\{\al\}}$ and using Fubini, we obtain 

$$A\le C\int_{(2Q)^c}|\partial^{t}(\varphi_Q(k^1*\psi))*\frac 1{|x|^{n+1-\{\al\}}})(x)|dx,$$
for some multiindex $t=(t_1,\cdots,t_n)$ with $|t|=n-[\al]-1$.

\noindent We will now use Leibniz formula \eqref{Leibniz} (with $s$ replaced by $t$) and the fact that

\begin{equation}\label{cota}
 |\partial^{t-r}k^1(x)| \le C\, |x|^{-(\alpha+|t|-|r|)}\,.
\end{equation}
Therefore, since  $\al+|t|-|r|<n$, by Fubini and H\"older
we obtain

\begin{equation}\label{lu}
\begin{split}
\int_{Q}|\partial^t(\varphi_Q(k^1*\psi))(y)|dy&\leq
 C\sum_{|r|=0}^{|t|}\int_{Q}|\partial^r\varphi_Q(y)||\partial^{t-r}(k^1*\psi)(y)|dy\\\\&\leq
 C\sum_{|r|=0}^{|t|}\int_{2Q}|\psi(x)|\int_Q\frac{|\partial^r\varphi_Q(y)|}{|y-x|^{\al+|t|-|r|}}dydx\\\\
&\leq C\|\psi\|_q l(Q)^{\frac n p+1-\{\al\}}.
\end{split}
\end{equation}

\noindent Hence, if we apply Fubini again, we get

$$A\le C \int_{Q}|\partial^t(\varphi_Q(k^1*\psi))(y)|\int_{(2Q)^c}\frac{dxdy}{|y-x|^{n+1-\{\al\}}}
\leq Cl(Q)^{\frac n p}\, \|\psi\|_q.$$

To estimate term $B$, we will use (\ref{Leibniz}) and for each $0\le |r|\le |t|$ we will add and substract 
$\partial^r\varphi_Q(x)\partial^{t-r}(k^1*\psi)(y)$ in the integral to gain integrability, namely

\begin{equation*}
\begin{split}
B&\le C\sum_{|r|= 0}^{|t|}\int_{2Q}|\int_Q\frac{(\partial^r\varphi_Q(y)-\partial^r\varphi_Q(x))\partial^{t-r}(k^1*\psi)(y)}{|y-x|^{n+1-\{\al\}}}dy|dx\\\\&+%C\sum_{|r|=0}^{|s|}
%\int_{2Q}|\partial^r\varphi_Q(x)||\int_Q\frac{\partial^{s-r}(k^1*\psi)(y)}{|y-x|^{n+1-\{\al\}}}dy|dx\\\\&=&B_1+B_2
C\sum_{|r|=0}^{|t|}\int_Q|\partial^r\varphi_Q(x)||\left(\Delta\partial^{t-r}k^1*\psi*\frac{1}{|y|^{n-1-\{\al\}}}\right)(x)|dx\\\\&
+C\sum_{|r|=0}^{|t|}\int_Q|\partial^r\varphi_Q(x)||\int_{Q^c}\frac{\partial^{t-r}(k^1*\psi)(y)}{|y-x|^{n+1-\{\al\}}}dy|dx\\\\&=B_1+B_2+B_3,
\end{split}
\end{equation*}
the last identity being a definition for $B_1$, $B_2$ and $B_3$.\newline

\noindent Since arguing as in (\ref{lu}), (recall that
$|t|=n-[\al]-1$), we obtain
\begin{equation*}
\int_Q|\partial^{t-r}(k^1*\psi)(y)|dy\le\|\psi\|_ql(Q)^{\frac n
p+1-\{\al\}+|r|}, \end{equation*}
by the mean value theorem and
Fubini, we get that

\begin{equation}\label{b1}
\begin{split}
B_1&\leq\sum_{|r|=0}^{|t|}\frac{C}{l(Q)^{|r|+1}}\int_{Q}|\partial^{t-r}(k^1*\psi)(y)|\int_{2Q}\frac{dxdy}{|y-x|^{n-\{\al\}}}\\\\
&\leq C\sum_{|r|=0}^{|t|}l(Q)^{-|r|-1}\|\psi\|_ql(Q)^{\frac n
p+1-\{\al\}+|r|}l(Q)^{\{\al\}}\\\\&\leq Cl(Q)^{\frac n p}\,
\|\psi\|_q.
\end{split}
\end{equation}

%To deal with term $B_2$, we divide it into two terms

%\begin{equation*}
%\begin{split}
%B_2&\le \sum_{|r|=0}^{|s|}\int_Q|\partial^r\varphi_Q(x)||\left(\partial^{s-r}k^1*\psi*\frac{1}{|y|^{n+1-\{\al\}}}\right)(x)|dx\\\\&+\sum_{|r|=0}^{|s|}\int_Q|\partial^r\varphi_Q(x)||\int_{Q^c}\frac{\partial^{s-r}(k^1*\psi)(y)}{|y-x|^{n+1-\{\al\}}}dy|dx=B_{21}+B_{22}.
%\end{split}
%\end{equation*}

We deal now with term $B_2$.
By computing the Fourier transform of  the convolution $\Delta\partial^{t-r}k^1*\psi*\frac{1}{|y|^{n-1-\{\al\}}}$, one can see that for $|r|=0$,

$$\left(\Delta\partial^t k^1*\psi*\frac{1}{|y|^{n-1-\{\al\}}}\right)(x)=c\psi+ cS_0(\psi)(x),$$
\noindent where $c$ is a constant, and $S_0$ is a smooth homogeneous convolution Calder\'on-Zygmund operator. For $|r|\geq 1$,
we obtain

$$\left(\Delta\partial^{t-r} k^1*\psi*\frac{1}{|y|^{n-1-\{\al\}}}\right)(x)=cS_r(\psi)(x),$$
with $S_r$ a convolution operator with kernel of homogeneity $-(n-|r|)$. 
Since Calder\'on-Zygmund operators preserve $L^q(\Rn)$, $1<q<\infty$, using Young's inequality % $\|k*\psi\|_{L^q(Q)}\le C \|k\|_{L^1(Q)}\|\psi\|_q$
to estimate the $L^q(Q)$-norm of the convolution $S_r(\psi)$ and H\"older, we obtain

\begin{equation*}
\begin{split}
%B_{21}&=
B_2&=\sum_{|r|=0}^{|t|}\int_Q|\partial^r\varphi_Q(x)||S_r(\psi)(x)|dx\\\\&
\le C\|\varphi_Q\|_p\|\psi\|_q+C\sum_{|r|=1}^{|t|}\|\partial^r\varphi_Q\|_p\|S_r(\psi)\|_{L^q(Q)}\\\\&
\le Cl(Q)^{\frac n p}\|\psi\|_q+\|\psi\|_q\sum_{|r|= 1}^{|t|}l(Q)^{-|r|}l(Q)^{\frac n p}l(Q)^{|r|}\leq Cl(Q)^{\frac n p}\, \|\psi\|_q.
\end{split}
\end{equation*}

Now we are only left with term %$B_{22}$.
$B_3$. Since $\partial^r\varphi_Q$ is supported on $Q$, we can substract
$\partial^r\varphi_Q(y)=0$ for $y\in 3Q\setminus Q$, $|r|\geq 0$. Then,
\begin{equation*}
\begin{split}
%B_{22}&
B_3&=\int_Q|\int_{Q^c}\partial^r\varphi_Q(x)\frac{\partial^{t-r}(k^1*\psi)(y)}{|y-x|^{n+1-\{\al\}}}dy|dx\\\\
   &\le\sum_{|r|=0}^{|t|}\int_Q|\int_{3Q\setminus Q}\partial^{t-r}(k^1*\psi)(y)\frac{\partial^r\varphi_Q(x)-\partial\varphi_Q(y)}{|y-x|^{n+1-\{\al\}}}dy|dx\\\\
&+\sum_{|r|=0}^{|t|}\int_Q|\int_{(3Q)^c}\partial^r\varphi_Q(x)\frac{\partial^{t-r}(k^1*\psi)(y)}{|y-x|^{n+1-\{\al\}}}dy|dx=%B_{22}^1+B_{22}^2
B_{31}+B_{32}.
\end{split}
\end{equation*}

Using the mean value theorem and proceeding as in (\ref{b1}),
$$%B_{22}^1&
B_{31}\le\sum_{|r|=0}^{|t|}\frac{C}{l(Q)^{|r|+1}}\int_{Q}\int_{3Q}\frac{|\partial^{t-r}(k^1*\psi)(y)|}{|y-x|^{n-\{\al\}}}dydx\le
Cl(Q)^{\frac n p}\|\psi\|_q.$$

\noindent Notice that by (\ref{cota}), for $y\in (3Q)^c$, 
$$|\partial^{t-r}(k^1*\psi)(y)|\leq\int_{2Q}\frac{\psi(z)}{|z-y|^{\al+|t|-|r|}}dz\le C \|\psi\|_q\,l(Q)^{\frac n p}\,l(Q)^{-\al-|t|+|r|}.$$

\noindent Therefore,
\begin{equation*}
\begin{split}
%B_{22}^2&
B_{32}&\le\sum_{|r|=0}^{|t|}\frac{C}{l(Q)^{|r|}}\int_{Q}\int_{(3Q)^c}\frac{|\partial^{t-r}(k^1*\psi)(y)|}{|y-x|^{n+1-\{\al\}}}dydx\\\\&\le
c\|\psi\|_ql(Q)^{\frac n
p}\sum_{|r|=0}^{|t|}\frac{l(Q)^{-\al-|t|+|r|}}{l(Q)^{|r|}}\int_Q\int_{(3Q)^c}\frac
1{|y-x|^{n+1-\{\al\}}}dydx\\\\&\le Cl(Q)^{\frac n p}\|\psi\|_q,
\end{split}
\end{equation*}
which finishes the proof of (\ref{condition}) for the case
$\{\al\}>0.\qed$ 
\end{enumerate}

\section{A continuity property for the capacity $\ga$}

In this section we prove a continuity property for the capacity
$\ga$, $0<\alpha<n$, which will be used in the proof of Theorem
\ref{alfa}.

\begin{lemma}\label{extregalfa}
Let $\{E_j\}_j$ be a decreasing sequence of compact sets, with intersection the compact set $E\subset\Rn$ and let $0<\alpha<n$. Then $$\ga(E)=\lim_{j\to\infty}\ga(E_j).$$
\end{lemma}

%\begin{lemma}\label{extregalfai}
%Let $\{E_j\}_j$ be a decreasing sequence of compact sets, with intersection the compact set $E\subset\Rn$ and let $0<\alpha<n$. Then $$\gai(E)=\lim_{j\to\infty}\gai(E_j),$$ for each component $1\leq i\leq n.$
%\end{lemma}
{\em Proof.} Since, by definition,  the set function $\ga$ in non-decreasing,
$$\lim_{j\to\infty}\ga(E_j)\geq\ga(E)\,,$$
and the limit clearly exists. For each $j\geq 1$, let $T_j$ be a
distribution such that the potential $x/|x|^{1+\al}*T_j$ is in the
unit ball of $L^\infty(\Rn)$\,,  and

$$\ga(E_j)-\frac 1 j < |\langle T_j,1\rangle| \leq\ga(E_j).$$

\noindent We want to show that for each test function $\varphi$,
\begin{equation}\label{test}
\langle
T_j,\varphi\rangle\underset{j\to\infty}{\longrightarrow}\langle
T,\varphi\rangle,
\end{equation}
for some distribution $T$ whose potential $x/|x|^{1+\al}*T$ is in
the unit ball of $L^\infty(\Rn)$. If (\ref{test})
holds and $\varphi$ is a test function satisfying $\varphi\equiv
1$ in a neighbourhood of $E$, then
$$\lim_{j\to\infty}\ga(E_j)=\lim_{j\to\infty}|\langle T_j,1\rangle |=\lim_{j\to\infty}|\langle T_j,\varphi\rangle|= |\langle T,\varphi\rangle|\leq\ga(E).$$

To show (\ref{test}), set $k_\alpha(x)=1/|x|^{n-\al}$ and $f_j^i=x_i/|x|^{1+\al}*T_j$, $1\leq i\leq n$. We will treat first the case $n$ odd and of the 
form $n=2k+1$. By  (\ref{betterrepro})
\begin{eqnarray*}
\langle T_j,\varphi\rangle&=&c\sum_{i=1}^n \langle
f_j^i,\Delta^k\partial_i\varphi*k_\al\rangle\\\\&=&c\sum_{i=1}^n\int
f_j^i(x)\left(\Delta^k\partial_i\varphi*k_\al\right)(x)dx
\end{eqnarray*}

We mention here that if $n=2k$, then one argues in the same way, but one has to use another reproduction
formula analogous to (\ref{betterrepro}) for this case (see
\cite[Lemma 3.1]{laura1}).

Passing to a subsequence, we can assume that for each $1\leq i\leq n$, $f_j^i\longrightarrow
f^i$  in the weak $*$ topology of $L^\infty(\Rn)$. But then
$f_j^i(x)(\Delta^k\partial_i\varphi*k_\al)(x) \longrightarrow f^i(x)(\Delta^k\partial_i\varphi*k_\al)
(x)\,,\quad x \in \Rn \,.$ Since this pointwise convergence is bounded,
%because is $x\in$spt$\varphi$, then $$|(f_j^i(x)\Delta^k\partial_i\varphi*k_\al)
%(x)|\leq\\int\frac{|\Delta^k\partial_i\varphi(z)|}{|z-y|^{n-\al}}dz\|f_j^i\|_{\infty}\le Cl(Q)^\al, i si està a fora la x, doncs li carreguem derivades al nucli i ja està.
the dominated convergence theorem
yields

\begin{equation*}
\begin{split}
\lim_{j\to\infty}<T_j,\varphi>&=c\,\lim_{j\to\infty}\sum_{i=1}^n\int f_j^i(x)\left(\Delta^k\partial_i\varphi*k_\al\right)(x)dx
\\\\&=c\sum_{i=1}^n\int f^i(x)\left(\Delta^k\partial_i\varphi*k_\al\right)(x)dx.
\end{split}
\end{equation*}

Define the distribution $T$ by

$$<T,\varphi>=c\sum_{i=1}^n\int f^i(x)\left(\Delta^k\partial_i\varphi*k_\al\right)(x)dx.$$

Now we want to show that for $1\leq i\leq n$, $f^i=x_i/|x|^{1+\alpha}*T$. For that we regularize $f_j^i$
and $T_j$\,. Take $\chi \in{\cal C}_0^{\infty}(B(0,1))$ with $\int
\chi(x)\,dx=1$ and set $\chi_\ep(x)=\ep^{-n}\chi(x/\ep)$\,. Then
we have , as $j \rightarrow \infty$\,,

$$ \left(\chi_\ep*\frac{x_i}{|x|^{1+\alpha}}*T_j\right)(x)= \left(\chi_\ep*f_j^i\right)(x)
\longrightarrow \left(\chi_\ep*f^i\right)(x)\,,\quad x \in \Rn\,,$$
because $f_j^i$ converges to $f^i$ weak $*$ in $L^\infty(\Rn)\,.$ On
the other hand, since $\chi_\ep *\frac{x_i}{|x|^{1+\alpha}} \in \cc^\infty(\Rn)$ and
$T_j$ tends to $T$ in the weak topology of distributions, with
controlled supports, we have
$$
\left(\chi_\ep*\frac{x_i}{|x|^{1+\al}}*T_j\right)(x) \longrightarrow
\left(\chi_\ep*\frac{x_i}{|x|^{1+\al}}*T\right)(x)\,, \quad x \in \Rn\,.
$$
Hence
$$
\chi_\ep*\frac{x_i}{|x|^{1+\al}}*T = \chi_\ep*f^i\,, \quad \ep > 0\,,
$$
and so, letting $\ep \rightarrow 0$\,, $\frac{x_i}{|x|^{1+\al}}*T = f^i\,.$ \qed

\section{Sketch of the proof of Theorem \ref{alfa}.}

This section will be devoted to the proof of inequality
(\ref{mainineq}), namely
 $$\ga(E)\leq C\gam(E).$$

\noindent We will adapt the line of reasoning in \cite{semiad} and \cite{semiad2}, where Tolsa proves the semiadditivity of analytic capacity and 
continuous analytic capacity respectively. We will also use the modifications introduced in \cite{volberg}, where the semiadditivity of Lipschitz harmonic capacity is proven (see also \cite{tolsaaleix}).

In fact, when one analizes the proofs of \cite{semiad}, \cite{semiad2} and \cite{volberg} one realizes that they depend on two technical facts, 
the exterior regularity property of $\ga$ (see Lemma \ref{extregalfa}) and an $L^\infty-$localization result, 
which is Theorem \ref{localization1} in our setting. 
We must mention that the positivity properties of the symmetrization method for the Cauchy kernel discovered in \cite{me} and \cite{meve} are an 
important ingredient for the proofs of \cite{semiad} and \cite{semiad2}. 
In \cite{volberg} one has to circumvent this lack of positivity and modify Tolsa's idea.

We will explain now how each of the above mentioned main ingredients
take part in the proof of (\ref{mainineq}): As we proved in Lemma
\ref{extregalfa}, the capacities $\ga$, $0<\al<n$, enjoy the
exterior regularity property. This is also true for the capacities
$\gam$, $0<\al<n,$ just by the weak $\star$ compactness of the the
set of positive
measures having total variation not exceeding $1\,.$  %As in
%\eqref{GammaopGamma+}, one has
%\begin{equation}\label{comppositiu}
%C^{-1}\, \Gamma_{\hat k,\op}(E) \le \Gamma_{\hat k,+}(E) \le C\,
%\Gamma_{\hat k,\op}(E)\,.
%\end{equation}
We therefore can approximate a general compact set $E$ by sets which are
finite unions of cubes of the same side length in such a way
that the capacities $\ga$ and $\gam$ of
the approximating sets are as close as we wish to those of $E\,.$
Thus we can assume, without loss of generality,  that $E$ is a
finite union of cubes of the same size. %Since Tolsa's argument
%is based on an induction on scales one needs a starting point,
%which is the scale corresponding to the size of the squares
%forming $E\,.$
This will allow to implement an induction argument on the size of certain ($n$-dimensional) rectangles. The first step involves rectangles of diameter comparable to the side length of the cubes whose union is $E$.

The starting point of the general inductive step in the proof of
Tolsa's Theorem in~\cite{semiad} (and \cite{semiad2}) and in \cite{volberg} for the Lipschitz harmonic capacity case, consists in
the construction of a positive Radon measure $\mu$ supported on a
compact set $F$ which approximates $E$ in an appropriate sense.  The
construction of $F$ and $\mu$ gives readily that
\begin{equation}\label{mu}
\ga(E)
\le C\, \mu(F),
\end{equation}
and
\begin{equation}\label{mu2}
\gam(F) \le C\,\gam(E), \end{equation} which tells us that $F$ is
not too small but also not too big.
 However, one cannot expect, in the context of \cite{semiad} and
\cite{semiad2}  the Cauchy singular integral to be
 bounded on~$L^2(\mu)$. In our case one cannot expect the $\al-$Riesz
operator $R(\mu)$ to be bounded on~$L^2(\mu)$. One has to carefully
look for a compact subset $G$ of $F$ such that
\begin{itemize}
\item $\mu(F) \le C\,\mu(G)$.
\item The restriction $\mu_G$ of $\mu$ to $G$ has $\al-$growth.
\item The operator $R(\mu_G)$, is bounded on $L^2(\mu_G)$ with
dimensional constants.
\end{itemize}

\noindent Moreover, recall from (\ref{gaopgam}), that one has

$$C^{-1}\,\gaop(E)\le\gam(E)\le C\,\gaop(E).$$

\noindent This completes the proof because then
\begin{equation*}
\begin{split}
\ga(E) &\le C\, \mu(F) \le C\, \mu(G) \le
C\,\gaop(G)  \le C\,\gaop(F) \\*[5pt]
& \le C\,\gam(F) \le C\, \gam(E) \le
C\,\gaop(E) .
\end{split}
\end{equation*}

In \cite{semiad}, \cite{semiad2} and \cite{volberg} the set $F$ is defined as the
union of a special family of cubes $\{Q_i\}_{i=1}^N$ that cover the
set $E$ and approximate $E$ at an appropriate intermediate scale.
One then sets
$$F=\bigcup_{i=1}^NQ_i.$$

The construction of the set $F$ is different in the analytic
capacity case and in the Lipschitz harmonic capacity case. In
Tolsa's proof this construction is performed by using the positivity
properties of the symmetrization of the Cauchy kernel discovered in
\cite{me} and \cite{meve}. In our setting, the symmetrization of the
Riesz kernels $x/|x|^{1+\al}$ only gives a positive quantity for
$0<\alpha\leq 1$, (see \cite[Lemma 4.2]{laura1}), therefore we have
to circumvent the use of this positivity property and therefore
modify Tolsa's idea.  For this modification we will follow chapter 5
of \cite{volberg}, where this was done for the Lipschitz harmonic
capacity case, namely for $\alpha=n-1$. In fact the arguments in
chapter 5 of \cite{volberg} are written for more general
Calder\'on-Zygmund kernels of homogeneity $-\al$ and so they also work in
our setting. Therefore, the construction of the approximating set $F$
with properties (\ref{mu}) and (\ref{mu2}) can be done just as in
\cite[ch. 5]{volberg}.\newline
%\noindent Just by construction, the set $F$ satisfies $\gam(F)\leq C\,\gam(E)$.
%Following Tolsa's argument, to prove the comparability between the capacities $\ga$ and $\gam$ we must show that there exists some positive measure $\mu\in L_\al(E)$ with
%$\ga(E) \approx \mu(E)\,,$
%and such that, $R_\al(\mu)$ is bounded on $L^2(\mu)$ with absolute constants.
%Then $$\ga(E)\le C\,\mu(E)\le C\, \gaop(E)$$
% and by (\ref{comppositiu}) we are done.

%In order to construct $\mu$, we will need to apply an induction argument: We will prove that $\ga(E\cap R)\leq C\gam(E\cap R)$ for all $n-$dimensional rectangles $R$, by induction on the size of $R$.

% Our aim now is to construct a family of cubes $\{Q_i\}_{i=1}^N$ with bounded overlapping, satisfying the following properties:
%\begin{equation}\label{u}
%\gam(\bigcup_{i=1}^NQ_i)\leq C_0\gam(E),
%\end{equation}
%\begin{equation}\label{dos}
%\sum_{i=1}^N\gam(2Q_i\cap E)\leq C_1\gam(E),
%\end{equation}
%\begin{equation}\label{tres}
%\mbox{diam}(Q_i)\leq\frac 1{10}\mbox{diam}(E).
%\end{equation}
%The constants $C_0$ and $C_1$ are dimensional constants.
%Set $$F=\bigcup_{i=1}^NQ_i.$$

%Property (\ref{u}) tells us that $\gam(F)\leq C_0\gam(E)$. Tolsa's
%argument will allow us to show that for the approximating set $F$
%the following holds:

%\begin{equation}\label{mu2}
%\ga(E)\approx\mu(F)
%\end{equation}

%We will now illustrate how to obtain (\ref{mu2}).
By the definition of the capacity $\ga$  it follows
that there exists a real distribution $T_0$ supported on the compact
set $E$ such that
\begin{enumerate}
 \item $\displaystyle{|<T_0,1>|\geq \frac{\ga(E)}{2}.}$
 \item $T_0$ has $\al-$growth and $G_\alpha(T_0) \le 1$\,.
 \item $\displaystyle{\|\frac{x_j}{|x|^{1+\al}}* T_0\|_\infty\leq 1,}\,\,\,\,\,\,1\le j\le n$.
\end{enumerate}

The construction of $\mu$ is performed simultaneously with that of a
real measure $\nu$, which should be viewed as a model for $T_0$
adapted to the family of cubes $\{Q_i\}_{i=1}^N\,.$  For each cube
$Q_i$ choose a ball $B_i$ concentric with $Q_i$ with
radius $r_i$ comparable to $\ga(E\cap 2Q_i)$ and set
$$\mu=\sum_{i=1}^N\frac{r_i^\alpha}{{\cal L}^n(B_i)}{\cal L}^n_{|B_i}.$$

Consider now functions $\varphi_i\in{\cal C}_0^{\infty}(2Q_i)$,
$0\leq\varphi_i\leq 1$, $\|\partial^s \varphi_i\|_\infty\leq
C\,l(Q_i)^{-|s|}\,,$ and
$\sum_{i=1}^N\varphi_i=1$ on $\bigcup_iQ_i$. The measure $\nu$ is
defined as

$$\nu=\sum_{i=1}^N\frac{<T_0,\varphi_i>}{{\cal L}^n(B_i)}{\cal L}^n_{|B_i},$$
with ${\cal L}^n$ being the $n$-dimensional Lebesgue measure.

Notice that supp$(\nu)\subset\mbox{supp}(\mu)\subset F$. Moreover we
have $d\nu=bd\mu$, with
$\displaystyle{b=\frac{<T_0,\varphi_i>}{r_i^\al}}$ on $B_i$. At this
point, we need to show that our function $b$ is bounded to apply
later a suitable $T(b)$ Theorem. To estimate $\|b\|_{\infty}$ we
need the localization Lemma \ref{localization1}, proved in section
$3.2$, which gives us

$$\|\frac{x_j}{|x|^{1+\al}}*\varphi_i T_0\|_\infty\leq C\,,\,\,\,\,1\le j\le n\,. $$
We therefore obtain, by the definition of $\ga$,

\begin{equation}\label{capita}
 |<T_0 , \varphi_i>|\leq C\ga(2Q_i\cap E),\,\,\,\mbox{ for }\,1\le i\le N\,.
\end{equation}
Hence $\|b\|_\infty\leq C$.  It is now easy to see why $\ga(E)\leq C\mu(F)$:
$$\ga(E)\leq 2|\langle T_0,1\rangle|\le 2\sum_{i=1}^N |\langle
T_0,\varphi_i\rangle|\leq C\sum_{i=1}^N\ga(E\cap
2Q_i)=C\mu(F).$$

Now everything is ready to apply a suitable
variant of the $T(b)$ Theorem (see \cite{ntv}). There is still one more difficulty, in applying the Nazarov, Treil and Volberg $T(b)-$type theorem, one needs finding a substitute for what they call the suppressed operators. It was already explained in \cite{laura1} that there are at least two versions of such operators for the Riesz kernels that work appropriately (see \cite[(2.7)  and (2.13)]{laura1}).

\noindent
Departament de Matem\`atiques, Universitat Aut\`onoma de
Barcelona, 08193 Bellaterra (Barcelona), Catalunya.\newline\newline
{\em E-mail:} {\tt laurapb@mat.uab.cat}

\end{document}